\documentclass[reqno]{alea2}
\usepackage{natbib}
\usepackage{fancyhdr}
\usepackage{graphicx}

\eheader{Alea}{{\bf }}{0}{0}{0}
\usepackage{amssymb}
\usepackage{amsmath}

\renewcommand{\cite}{\citet}

\makeatletter \@addtoreset{equation}{section} \makeatother

\renewcommand\theequation{\thesection.\arabic{equation}}
\renewcommand\thefigure{\thesection.\@arabic\c@figure}
\renewcommand\thetable{\thesection.\@arabic\c@table}

\newtheorem{theorem}{Theorem}[section]
\newtheorem{lemma}[theorem]{Lemma}
\newtheorem{proposition}[theorem]{Proposition}
\newtheorem{corollary}[theorem]{Corollary}

\newtheorem{definition}[theorem]{Definition}
\newtheorem{condition}[theorem]{Condition}

\newtheorem{remark}[theorem]{Remark}
\newtheorem{example}[theorem]{Example}

\def \hat{\widehat}
\def \bar{\overline}
\def \reff#1{(\ref{#1})} 
 
\def \R{{\mathbb R}}
\def \Z{{\mathbb Z}} 
 
\def \one{{\bf 1}\hskip-.5mm}

\begin{document}

\date{, accepted }
\keywords{Poisson point processes, spatial 
point processes, birth and death processes, Poisson 
random measures, stochastic  
equations, ergodicity. } 
\subjclass{Primary:  60K35, 60G55 Secondary:  60J27, 60J35, 60H20,
  82B21, 82C21. } 

\author{Nancy L. Garcia and Thomas G.  Kurtz}
\address{ Departamento de 
Estat\'\i stica, IMECC -- UNICAMP, Caixa Postal 6065, 13.081-970 -
Campinas, SP -  BRAZIL.}
\email{nancy@ime.unicamp.br} 
\address {Department of 
Mathematics, University of Wisconsin, Madison 
Wisconsin 53706, U.S.A.}
\email{kurtz@math.wisc.edu}

\title[] {Spatial birth and death processes as solutions of 
stochastic equations} 

\begin{abstract}
  Spatial birth and death processes are obtained as 
solutions of a system of stochastic equations.  The 
processes are required to be locally finite, but may 
involve an infinite population over the full (noncompact) 
type space.  Conditions are given for existence and 
uniqueness of such solutions, and for temporal and 
spatial ergodicity.  For birth and death processes with 
constant death rate, a sub-criticality condition on the 
birth rate implies that the process is ergodic and 
converges exponentially fast to the stationary 
distribution.  

\end{abstract}

\maketitle

\setcounter{equation}{0}

\section{Introduction} 

Spatial birth and death processes in which the birth and 
death rates depend on the configuration of the system 
were first studied by Preston (1975).  His approach was 
to consider the solution of the backward Kolmogorov 
equation, and he worked under the restriction that there 
were only a finite number of individuals alive at any 
time.  Under certain conditions, the processes exist and 
are temporally ergodic, that is, there exists a unique 
stationary distribution.  The more general setting 
considered here requires only that the number of points 
alive in any compact set remains finite at all times.  

Specifically, we assume that our population is 
represented as a countable subset of points in 
a complete, separable 
metric space $S$  (typically, $S\subset\R^d$). 
We will identify the 
subset with the counting measure $\eta$ given by 
assigning unit mass to each point, that is, $\eta (B)$ is the 
number of points in a set $B\in {\mathcal B}(S)$.  
(${\mathcal B}(S)$ will denote the Borel subsets of 
$S$.)
We will use the terms 
point process and random counting measure 
interchangeably.  
With this identification in mind, let 
${\mathcal N}(S)$ be the collection of counting measures on the
metric space $S$.  The state space for 
our process will be some subset of ${\mathcal N}(S)$.   All processes 
and random variables are defined on a complete probability space 
$(\Omega ,{\mathcal F},P)$.  

The spatial birth and death process is specified in terms 
of non-negative functions $\lambda :S\times {\mathcal N}(S)\rightarrow 
[0,\infty )$ and 
$\delta :S\times {\mathcal N}(S)\rightarrow [0,\infty )$ and a reference measure $
\beta$ on $S$ 
(typically Lebesgue measure $m_d$, if $S\subset {\Bbb R}^d$).  $
\lambda$ is the 
birth rate and $\delta$ the death rate.  If the point 
configuration at time $t$ is $\eta\in {\mathcal N}(S)$, then the probability 
that a point in a set $B\subset S$ is added to the configuration 
in the next time interval of length $\Delta t$ is approximately 
$\int_B\lambda (x,\eta )\beta (dx)\Delta t$ and the probability that a point $
x\in\eta$ is 
deleted from the configuration in the next time interval 
of length $\Delta t$ is approximately $\delta (x,\eta )\Delta t$.  Under these 
assumptions, the generator of the process should be of 
the form 
\begin{equation}AF(\eta )=\int (F(\eta +\delta_x)-F(\eta ))\lambda 
(x,\eta )\beta (dx)+\int (F(\eta -\delta_x)-F(\eta ))\delta (x,\eta 
)\eta (dx)\label{gener}\end{equation}
for $F$ in an appropriate domain.  

Following the work of Preston, spatial birth and death 
processes quickly found application in statistics when 
Ripley (1977) observed that spatial point patterns could 
be simulated by constructing a spatial birth and death 
process having the distribution of the desired pattern as 
its stationary distribution and then simulating the birth 
and death process for a long time, a procedure now 
known as Markov chain Monte Carlo. 

The two best-known classes of spatial point processes 
are Poisson random measures and Gibbs distributions.
 
\subsection{Poisson random measures}  
Let $\beta$ be a $\sigma$-finite measure on $S$, $(S,d_S)$ a complete, 
separable metric space.  $\xi$ is 
a Poisson random measure on $S$ with mean 
measure $\beta$ if for each $B\in {\mathcal B}(S)$, $\xi (B)$ has a 
Poisson distribution with expectation $\beta (B)$ and $\xi (B)$ and 
$\xi (C)$ are independent if $B\cap C=\emptyset$.  Taking $\lambda 
=\delta\equiv 1$,
then the Poisson random measure with mean measure $\beta$ 
gives the unique stationary distribution for the birth and 
death process with generator
\begin{equation}AF(\eta )=\int (F(\eta +\delta_x)-F(\eta ))\beta 
(dx)+\int (F(\eta -\delta_x)-F(\eta ))\eta (dx).\label{gen2}\end{equation}
  Letting $\mu_{\beta}^0$ denote this distribution, the 
stationarity can be checked by verifying that
\[\int_{{\mathcal N}(S)}AF(\eta )\mu_{\beta}^0(d\eta )=0.\]
This assertion follows from the standard identity
\begin{equation}E[\int_Sh(\xi -\delta_x,x)\xi (dx)]=E[\int_Sh(\xi 
,x)\beta (dx)].\label{poisid}\end{equation}
See Daley and Vere-Jones (1988), p. 188, Equation (6.4.11).

\subsection{Gibbs distributions} Assume that $\beta (S)<\infty$.  
Consider the class of spatial point processes specified 
through a density (Radon-Nikodym derivative) with 
respect to a Poisson point process with mean measure $\beta$, 
that is, the distribution of the point process is given by 
\begin{equation}\mu_{\beta ,H}(d\eta )=\frac 1{Z_{\beta ,H}}e^{-H
(\eta )}\mu_{\beta}^0(d\eta ),\label{ccol1}\end{equation}
where $H(\eta )$ is referred to as the {\em energy function}, $Z_{
\beta ,H}$ is a 
normalizing constant, and $\mu_{\beta}^0$ is the law of a Poisson 
process with mean measure $\beta$.  Therefore, the 
state space for this process is ${\mathcal S}=\{\eta\in {\mathcal N}(S);H
(\eta )<\infty \}$, 
the set of configurations with positive density.  We 
assume that $H$ is hereditary in the sense of Ripley 
(1977),  that is $H(\eta )<\infty$ and $\tilde{\eta}\subset\eta$ implies $
H(\tilde{\eta })<\infty$.  
Ripley showed that such a measure $\mu_{\beta ,H}$ is the 
stationary
distribution of a spatial birth and death process.  In fact, 
there is more than one birth and death
process that has $\mu_{\beta ,H}$ as a
stationary distribution; we simply require that $\lambda (x,\eta 
)>0$ if 
$H(\eta +\delta_x)<\infty$ and that 
$\lambda$ and $\delta$ satisfy 
\begin{equation}\lambda (x,\eta )e^{-H(\eta )}=\delta (x,\eta +\delta_
x)e^{-H(\eta +\delta_x)}.\label{eqcol1}\end{equation}
This equation is a detailed balance condition which 
ensures that births from $\eta$ to $\eta +\delta_x$ match deaths from 
$\eta +\delta_x$ to $\eta$ and that the process is time-reversible with 
\reff{ccol1} as its stationary distribution.  
Again, this assertion can be verified by showing that
\[\int AF(\eta )\mu_{\beta ,H}(d\eta )=\frac 1{Z_{\beta ,H}}\int 
AF(\eta )e^{-H(\eta )}\mu^0_{\beta}(d\eta )=0.\]
This 
identity again follows from (\ref{poisid}).

Notice that equation (\ref{eqcol1}) says that any pair of 
birth and death rates such that 
\[\frac {\lambda (x,\eta )}{\delta (x,\eta +\delta_x)}=\exp\{-H(\eta 
+\delta_x)+H(\eta )\}\]
will give rise to a process with stationary distribution
given by \reff{ccol1}.  We can always 
take $\delta (x,\eta )=1$, that is, whenever a point is added to the 
configuration, it lives an exponential length of time 
independently of the configuration of the process.  

For example, consider a spatial point process on a 
compact set $S\subset\R^d$ given by a Gibbs distribution with 
pairwise interaction potential $\rho (x_1,x_2)\ge 0$, that is, for 
$\eta =\sum_{i=1}^m\delta_{x_i}$,
\begin{eqnarray}
H_{\rho}(\eta )&=&\sum_{i<j}\rho (x_i,x_j)\label{eqcol1a}\\
&=&\frac 12[\int\int\rho (x,y)\eta (dx)\eta (dy)-\int\rho (x,x)\eta 
(dx)]\}\nonumber\end{eqnarray}
and the distribution of the point process is absolutely 
continuous with respect to the spatially homogeneous 
Poisson process with constant intensity $1$ (or equivalently 
with Lebesgue mean measure) on $S$.  
Taking $\delta (x,\eta )\equiv 1$ and $\lambda (x,\eta )=\exp\{-\int
\rho (x,y)\eta (dy)\}$, 
the distribution determined by 
(\ref{eqcol1a}) is the stationary distribution for the birth 
and death process 
with infinitesimal generator 
\begin{equation}\label{eqcol2}AF(\eta )=\int e^{-\int\rho (x,y)\eta 
(dy)}(F(\eta +\delta_x)-F(\eta ))dx+\int (F(\eta -\delta_x)-F(\eta 
))\eta (dx).\end{equation}

Another 
example is the area-interaction point process introduced by 
Baddeley and Van Lieshout (1995).  This point process is 
absolutely continuous with respect to the spatial Poisson 
process with Lebesgue mean measure $m_d$ on $S\subset\R^d$ and
 $H(\eta )=\eta (S)\log\rho -m_d(\eta\oplus G)$, so the
Radon-Nikodym derivative is given by
\begin{equation}L(\eta )=\frac 1Z\rho^{\eta (S)}\gamma^{-m_d(\eta
\oplus G)}.\label{eqcol3}\end{equation}
Again, $Z$ is the normalizing constant, $\rho$ and $\gamma$ are 
positive parameters, and $G$ is a compact (typically 
convex) subset of $\R^d$ referred to as the {\em grain}.  The set 
$\eta\oplus G$ is given by 
\[\eta\oplus G=\cup \{x\oplus G;x\in\eta \}.\]
The parameter $\gamma$ controls the area-interaction among 
the points of $\eta$:  the process is \emph{attractive} if 
$\gamma >1$ and \emph{repulsive} otherwise.  (See Lemma 
\ref{mono}.)  If $\gamma =1$ the 
point process is just the Poisson random measure with mean 
measure  $\rho m_d$.  The case $\gamma >1$ is related 
to the \emph{Widow-Rowlinson model} introduced by 
Widow and Rowlinson (1970).  The case of 
\emph{area-exclusion} corresponds to a suitable limit 
$\gamma\to 0$.  A birth and death process with
stationary distribution given by the area-interaction 
distribution can be obtained by taking the 
unit death rate and the birth rate given by 
\begin{equation}\label{eqcol4}\lambda (x,\eta )=\rho\,\gamma^{-m_d(
(x+G)\setminus (\eta\oplus G))}.\end{equation}

\subsection{Overview} 
The spatial birth and death processes that correspond to 
the Gibbs distributions discussed above involve finite 
configurations and it is straightforward to see that they 
are uniquely characterized by their birth and death 
rates, for example, as solutions of the martingale 
problem associated with the generator $A$ given in 
(\ref{gener});
however, if the configurations are infinite and the 
total birth and death rates are infinite,
the existence and uniqueness of the processes are not so 
clear. 
In Section \ref{section2}, we represent these processes 
as solutions of a system of stochastic equations and give 
conditions for existence and uniqueness of solutions for 
the equations as well as for the corresponding 
martingale problems.  These equations are very 
useful in studying the asymptotic behavior of the birth and 
death processes,
including temporal and/or spatial-ergodicity and the speed of 
convergence to the stationary distribution.  

The uniqueness conditions given here are direct analogs 
of Liggett's (1972) conditions for existence and uniqueness 
for lattice indexed interacting particle systems.  
Stochastic equations for lattice indexed systems were 
formulated in Kurtz (1980) using time-changed Poisson 
processes and existence and uniqueness given under 
Liggett's conditions.  Stochastic equations for spatial 
birth and death processes of the type considered here 
were formulated in Garcia (1995) using a spatial version 
of the time-change approach.  Existence and 
uniqueness were again given under analogs of Liggett's 
conditions.  

One disadvantage to the time-change approach taken in 
Kurtz (1980) and Garcia (1995) is that the filtration to 
which the process is adapted depends on the solution.  
A stochastic equation that 
avoids this difficulty can be formulated by representing 
the birth process as a thinning of a Poisson random 
measure.  Intuitively, the approach is analogous to the 
rejection method for simulating random variables.  The 
fact that counting processes and more general marked 
counting processes can be obtained by thinning Poisson 
random measures is well known, particularly in the 
context of simulation.  (See, for example, Daley and 
Vere-Jones (2003), Section 7.5.) Stochastic equations exploiting 
this approach were formulated for 
lattice systems in Kurtz and Protter (1996) and for 
general spatial birth processes by Massouli\'e (1998).  In 
both cases, uniqueness was obtained under conditions 
analogous to Liggett's.

Section \ref{sectiondeath} considers temporal 
ergodicity for birth and death 
processes in noncompact $S$ (more precisely, $S$ and $\beta$ 
with $\beta (S)=\infty$) and spatial ergodicity 
for $S=\R^d$ and translation invariant 
birth rates. 
We give conditions for  ergodicity and 
exponential convergence to the stationary distribution.
It is well known that these processes are 
temporally and spatially ergodic if the birth and death 
rates are constant (the stationary measure being 
Poisson).  More generally, in Theorem \ref{thte}, we show 
that 
if the birth rate  
satisfies the conditions of Theorem \ref{thmeu} with 
$M<1$ and the death rate is constant ($\delta\equiv 1$),
the system is temporally ergodic and for every initial 
distribution, the distribution of the solution converges at an 
exponential rate to the stationary distribution.  
For $S=\R^d$ and $\lambda$ translation invariant, from the 
stochastic equation, we see that spatial ergodicity of the 
initial distribution implies spatial ergodicity of the 
solution at each time $0<t<\infty$.
Unfortunately, it is not 
clear, in general, how to carry this 
conclusion over to $t=\infty$, that 
is, to the limiting distribution of the solution, although 
 in the case $M<1$, 
spatial ergodicity holds for  the unique stationary 
distribution as well.
 We give some additional
conditions under which spatial ergodicity of 
the limiting distribution can be obtained.

Fern\'{a}ndez, 
Ferrari and Garcia (2002) study ergodicity of spatial birth 
and death processes using 
a graphical representation 
to construct the stationary distribution that is closely 
related to the stochastic equations we consider here. 
They give conditions for an exponential rate of 
convergence to the stationary distribution and for 
spatial ergodicity of the stationary distribution similar 
to those given here, but for a more restricted class of 
models.

Throughout, $\bar {C}(S)$ will denote the space of bounded 
continuous functions on $S$ and ${\mathcal B}(S)$ the Borel subsets of 
$S$.

The stochastic equations we consider will be driven by a Poisson
random measure $N$ on $U\times [0,\infty )$ for an appropriate space
$U$, having mean measure of the form $\nu\times m_1$, where $m_1$ is
Lebesgue measure on $[0,\infty )$.  Then for $B\in {\mathcal B} (U)$ with
$\nu (B)<\infty$, $N(B,t)\equiv N(B\times [0,t])$ is just an ordinary
Poisson process with intensity $\nu (B)$.  For a filtration $\{{\mathcal
  F}_t \}$, we say that $N$ is {\em compatible\/} with $\{{\mathcal
  F}_t\}$ if and only if for each $B\in {\mathcal B}(U)$ with $\nu
(B)<\infty$, $N(B,\cdot )$ is $\{{\mathcal F}_t\}$-adapted and
$N(B,t+s)-N(B,t)$ is independent of ${\mathcal F}_t$ for $s,t\geq 0$.

\setcounter{equation}{0}

\section{Spatial birth and death processes as solutions of 
stochastic  equations} \label{section2} 

A birth and death process as described in the previous 
section can be represented as the solution of a system 
of stochastic equations.  The approach is similar to 
Garcia (1995) where such processes were obtained as 
solutions of time-change equations.  We assume that the 
individuals in the birth and death process are 
represented by points in a Polish space $S$.  Typically, $S$ 
will be $\R^d$, $\Z^d$, or a subset of one of these, but we do 
not rule out more general spaces.  Let $K_1\subset K_2\subset\cdots$ 
satisfy $\cup_kK_k=S$, and let $c_k\in\bar {C}(S)$ satisfy $c_k\geq 
0$ and 
$\inf_{x\in K_k}c_k(x)>0$.  ${\mathcal N}(S)$ will denote the collection of 
counting measures on $S$ and ${\mathcal S}$ will denote 
$\{\zeta\in {\mathcal N}(S):\int_Sc_k(x)\zeta (dx)<\infty ,k=1,2,\ldots 
\}$.  Without loss of 
generality, we can assume that $c_1\leq c_2\leq\cdots$.  Let 
${\mathcal C}=\{f\in\bar {C}(S):|f|\leq ac_k\mbox{\rm \ for  some }k\mbox{\rm \ and }
a>0\}$, and topologize ${\mathcal S}$ 
by the weak$*$ topology generated by ${\mathcal C}$, that is, $\zeta_
n\rightarrow\zeta$ if 
and only if $\int_Sfd\zeta_n\rightarrow\int_Sfd\zeta$ for all $f\in 
{\mathcal C}$.  (Note that ${\mathcal C}$ is 
linear and that, with this topology, ${\mathcal S}$ is Polish.)  $D_{
{\mathcal S}}[0,\infty )$ 
will denote the space of cadlag ${\mathcal S}$-valued functions with 
the Skorohod ($J_1$) topology.  We assume that $\lambda$ and $\delta$ are 
nonnegative,
Borel measurable functions on $S\times {\mathcal C}$.

Let $\beta$ be a $\sigma$-finite, Borel measure on $S$.
   We assume 

\begin{condition}\label{bound}For 
each compact ${\mathcal K}\subset {\mathcal S}$, the birth rate $\lambda$ satisfies 
\begin{equation}\sup_{\zeta\in {\mathcal K}}\int_Sc_k(x)\lambda (x,\zeta 
)\beta (dx)<\infty ,\quad t>0,\quad k=1,2,\ldots ,\label{brest}\end{equation}
and 
\begin{equation}\delta (x,\zeta )<\infty ,\quad\zeta\in {\mathcal S},
\quad x\in\zeta .\label{ficond}\end{equation}
\end{condition}

We also assume that $\lambda$ and $\delta$ satisfy the following 
continuity condition.  

\begin{condition}\label{contcnd} 
If 
\begin{equation}\lim_{n\rightarrow\infty}\int_Sc_k(x)|\zeta_n-\zeta 
|(dx)=0,\label{ptwise}\end{equation}
for each $k=1,2,\ldots$, then 
\begin{equation}\lambda (x,\zeta )=\lim_{n\rightarrow\infty}\lambda 
(x,\zeta_n),\quad\delta (x,\zeta )=\lim_{n\rightarrow\infty}\delta 
(x,\zeta_n).\label{ccond}\end{equation}
\end{condition}

Note that since \reff{ptwise} implies $\zeta_n$ converges to $\zeta$ 
in ${\mathcal S}$, the continuity condition \reff{ccond} is weaker than 
continuity in ${\mathcal S}$; however, we have the following 
condition under which convergence in ${\mathcal S}$ implies 
(\ref{ptwise}).

\begin{lemma}\label{varconv}
Suppose $\zeta_0,\zeta_1,\zeta_2,\ldots\in {\mathcal S}$  and $\zeta_
n\leq\zeta_0$, $n=1,2,\ldots$.  If 
$\zeta_n\rightarrow\zeta$ in ${\mathcal S}$, then (\ref{ptwise}) holds.
\end{lemma}

\begin{proof}
$\zeta_n\leq\zeta_0$ implies that, considered as a measure, $\zeta_
n<<\zeta_0$ 
and $\frac {d\zeta_n}{d\zeta_0}=\frac {\zeta_n(\{x\})}{\zeta_0(\{
x\})}\leq 1$, almost everywhere $\zeta_0$.  Furthermore, 
$\zeta_n\rightarrow\zeta$ in ${\mathcal S}$ implies $\zeta_n(\{x\})\rightarrow
\zeta (\{x\})$ for each $x\in\zeta_0$, since the 
support of $\zeta_0$ consists of a countable collection of 
isolated points. Consequently, 
\[c_k(x)\geq c_k(x)|\frac {\zeta_n(\{x\})}{\zeta_0(\{x\})}-\frac {
\zeta (\{x\})}{\zeta_0(\{x\})}|\rightarrow 0,\]
and since $\int_Sc_k(x)\zeta_0(dx)<\infty$, the dominated convergence 
theorem implies
\[\lim_{n\rightarrow\infty}\int_Sc_k(x)|\zeta_n-\zeta |(dx)=\lim_{
n\rightarrow\infty}\int_Sc_k(x)|\frac {\zeta_n(\{x\})}{\zeta_0(\{
x\})}-\frac {\zeta (\{x\})}{\zeta_0(\{x\})}|\zeta_0(dx)=0.\]
\end{proof}

\begin{lemma}
Suppose ${\mathcal H}\subset {\mathcal S}$ and $\zeta_0\in {\mathcal S}$ satisfy $
\zeta\leq\zeta_0$, $\zeta\in {\mathcal H}$.  If $\delta$ 
satisfies (\ref{ficond})  and Condition \ref{contcnd}, then 
\[\sup_{\zeta\in {\mathcal H}}\delta (x,\zeta )<\infty .\]
\end{lemma}

\begin{proof}
${\mathcal H}$ is relatively compact in ${\mathcal S}$, so any sequence $
\{\zeta_n\}\subset {\mathcal H}$ 
has a subsequence that converges in ${\mathcal S}$ 
and, by 
Lemma \ref{varconv},
satisfies (\ref{ptwise}).  
 Fix $x\in S$, and 
let $\{\zeta_n\}$ satisfy $\lim_{n\rightarrow\infty}\delta (x,\zeta_
n)=\sup_{\zeta\in {\mathcal H}}\delta (x,\zeta )$. Then there 
is a subsequence that converges to some $\hat{\zeta}\in {\mathcal S}$ and 
hence, $\sup_{\zeta\in {\mathcal H}}\delta (x,\zeta )=\delta (x,\hat{
\zeta })<\infty$.
\end{proof}

\begin{lemma}
Suppose that for each $x\in S$, there exists $k(x)$ such that
 $\lambda (x,\zeta +\delta_y)=\lambda (x,\zeta )$ for $y\notin K_{
k(x)}$.  Then $\lambda$ satisfies 
Condition \ref{contcnd}\ and similarly for $\delta$.
\end{lemma}

\begin{proof} Note that $\zeta\in {\mathcal S}$ implies $\zeta (K_{k(x)})<\infty$  and that 
$\zeta_n\rightarrow\zeta$ implies that for $n$ sufficiently large, $
\zeta_n$ restricted 
to $K_{k(x)}$ coincides with $\zeta$ restricted to $K_{k(x)}$, so 
$\lambda (x,\zeta_n)=\lambda (x,\zeta )$. \hfill 
\end{proof}

Let $N$ be a 
Poisson random measure on $S\times [0,\infty )^3$ with mean 
measure $\beta (dx)\times ds\times e^{-r}dr\times du$.
 Let $\eta_0$ be an ${\mathcal S}$-valued random variable independent of $
N$, 
and let $\hat{\eta}_0$ be the point process on $S\times [0,\infty 
)$ obtained by 
associating to each ``count'' in $\eta_0$ an independent, unit 
exponential random variable, that is, for $\eta_0=\sum_{i=1}^{\infty}
\delta_{x_i}$, 
set 
\begin{equation}\hat{\eta}_0=\sum_{i=1}^{\infty}\delta_{(x_i,\tau_
i)},\label{etahat}\end{equation}
where the $\{\tau_i\}$ are independent unit exponentials, 
independent of $\eta_0$ and $N$.  The birth and death process $\eta$ 
should satisfy a stochastic equation of the form 
\begin{eqnarray}
\label{stoch1a}\eta_t(B)&=&\int_{B\times [0,t]\times [0,\infty )^
2}\one_{[0,\lambda (x,\eta_{s-})]}(u)\one_{(\int_s^t\delta (x,\eta_
v)\,dv,\infty )}(r)N(dx,ds,dr,du) \nonumber \\
&&\quad\quad\quad\quad\quad\quad\quad\quad +\int_{B\times [0,\infty 
)}\one_{(\int_0^t\delta (x,\eta_s)\,ds,\infty )}(r)\hat{\eta}_0(d
x,dr).\end{eqnarray}
To be precise, let $\eta$ be a process with sample paths in 
$D_{{\mathcal S}}[0,\infty )$ that is adapted to a filtration $\{{\mathcal F}_
t\}$ with respect 
to which $N$ is compatible.  (Note that \reff{brest} 
ensures that the integral with respect to $N$ on the right 
exists and determines an ${\mathcal S}$-valued random variable, and 
the continuity condition \reff{ccond} and the finiteness 
of $\delta (x,\zeta )$ ensure that $\delta (s,\eta_t)$ is a cadlag function of $
t$, so 
that the integrals $\int_s^t\delta (x,\eta_v)dv$ exist.)  Then $\eta$ is a 
solution of \reff{stoch1a} if and only if the identity 
\reff{stoch1a} holds almost surely for all $B\in {\mathcal B}(S)$ and 
$t\geq 0$ (allowing $\infty =\infty$).

\begin{lemma}
Suppose Condition \ref{bound} holds.  If $\eta$ is a solution of 
(\ref{stoch1a}), then for each $T>0$,
\begin{equation}\int_0^T\int_Sc_k(x)\lambda (x,\eta_s)\beta (dx)d
s<\infty\quad a.s.,\label{intbnd}\end{equation}
 $\eta^{*}_T$  defined by
\[\eta^{*}_T(B)=\int_{B\times [0,T]\times [0,\infty )^2}\one_{[0,
\lambda (x,\eta_{s-})]}(u)N(dx,ds,dr,du)\]
is an element of ${\mathcal S}$, 
\[\eta_t\leq\eta_T^{*}+\eta_0,\quad 0\leq t\leq T,\]
and 
\[\lim_{s\rightarrow t+}\int_Sc_k(x)|\eta_s-\eta_t|(dx)=0,\quad t
\geq 0.\]
\end{lemma}

\begin{proof}
Since for almost every $\omega\in\Omega$, the closure of 
$\{\eta_s:0\leq s\leq T\}$ is compact, Condition \ref{bound}\ implies
\[\sup_{s\leq T}\int_Sc_k(x)\lambda (x,\eta_s)\beta (dx)<\infty\quad 
a.s.,\]
and hence (\ref{intbnd}).  Letting 
\[\tau_c=\inf\{t:\int_0^t\int_Sc_k(x)\lambda (x,\eta_s)\beta (dx)
ds>c\},\]
we have
\[E[\int_Sc_k(x)\eta^{*}_{T\wedge\tau_c}(dx)]=E[\int_0^{T\wedge\tau_
c}\int_Sc_k(x)\lambda (x,\eta_s)\beta (dx)ds]\leq c,\]
and since $\lim_{c\rightarrow\infty}\tau_c=\infty$ a.s., it follows that 
$\int_Sc_k(x)\eta^{*}_T(dx)<\infty$, a.s.  implying $\eta^{*}_T\in 
{\mathcal S}$ a.s.  The last 
statement then follows by Lemma \ref{varconv}.
\end{proof}

If $\eta$ is a solution of \reff{stoch1a} and a point at $x$ was 
born at time $s\leq t$, then the ``residual clock time'' 
$r-\int_s^t\delta (x,\eta_v)dv$ is an ${\mathcal F}_t$-measurable random variable.  In 
particular, the counting-measure-valued process given by 
\begin{eqnarray}
\hat{\eta}_t(B\times D)&=&\int_{B\times [0,t]\times [0,\infty )^2}\one_{
[0,\lambda (x,\eta_{s-})]}(u)\one_D(r-\int_s^t\delta (x,\eta_v)\,
dv)N(dx,ds,dr,du)\nonumber\\
&&\quad\quad\quad\quad\quad\quad +\int_{B\times [0,\infty )}\one_
D(r-\int_0^t\delta (x,\eta_{s-})\,ds)\hat{\eta}_0(dx,dr)\label{stoch2a}\end{eqnarray}
is $\{{\mathcal F}_t\}$-adapted.  

Let $\hat {{\mathcal S}}$ denote the collection of counting measures $
\zeta$ on 
$S\times [0,\infty )$ such that $\zeta (\cdot\times [0,\infty ))\in 
{\mathcal S}$.  We can formulate an 
alternative equation for the $\hat {{\mathcal S}}$-valued process $\hat{
\eta}$ by 
requiring that 
\begin{eqnarray}
&&\int_{S\times [0,\infty )}f(x,r)\hat{\eta}_t(dx,dr)\label{stoch3a}\\
&&\qquad =\int_{S\times [0,\infty )}f(x,r)\hat{\eta}_0(dx,dr)\nonumber\\
&&\qquad\qquad\qquad +\int_{S\times [0,t]\times [0,\infty )^2}f(x
,r)\one_{[0,\lambda (x,\eta_{s-})]}(u)N(dx,ds,dr,du)\nonumber\\
&&\qquad\qquad\qquad -\int_0^t\int_{S\times [0,\infty )}\delta (x
,\eta_s)f_r(x,r)\hat{\eta}_s(dx,dr)ds,\nonumber\end{eqnarray}
for all $f\in\hat {{\mathcal C}}$, where $\hat {{\mathcal C}}$ is the collection of $
f\in\bar {C}(S\times [0,\infty ))$ 
such that $f_r\equiv\frac {\partial}{\partial r}f\in\bar {C}(S\times 
[0,\infty ))$, $f(x,0)=0$,   
$\sup_r|f(\cdot ,r)|,\sup_r|f_r(\cdot ,r)|\in {\mathcal C}$, and there exists $
r_f>0$ such 
that $f_r(x,r)=0$ for $r>r_f$.  Note that if $f\in\hat {{\mathcal C}}$ and
\begin{equation}f^{*}(x,r)=\int_0^r|f_r(x,u)|du,\label{absderiv}\end{equation}
then $f^{*}\in\hat {{\mathcal C}}$.
In \reff{stoch3a}, $\hat{\eta}_0$ can be 
any $\hat {{\mathcal S}}$-valued random variable that is independent of $
N$.

\subsection{Martingale problems} Let ${\mathcal D}(\hat {A})$ be the collection 
of functions $F$ of the form 
$F(\hat{\zeta })=e^{-\int_{S\times [0,\infty )}f(x,r)\hat{\zeta }
(dx,dr)}$, for non-negative  $f\in\hat {{\mathcal C}}$.  
Suppose that $\hat{\eta}$ is a solution of (\ref{stoch3a}) with 
sample paths in $D_{\hat {{\mathcal S}}}[0,\infty )$.  Assuming Condition \ref{bound}, 
\begin{equation}\int_0^t\int_Sc_k(x)\lambda (x,\eta_s)\beta (dx)d
s<\infty ,\quad k=1,2,\ldots .\label{intest}\end{equation}
By It\^o's formula 
\begin{eqnarray}
F(\hat{\eta}_t)\hspace{-2mm}&=&\hspace{-2mm}F(\hat{\eta}_0)+\int_{S\times [0,t]\times [0,\infty 
)^2}F(\hat{\eta}_{s-})(e^{-f(x,r)}-1)\one_{[0,\lambda (x,\eta_{s-}
)]}(u)N(dx,ds,dr,du)\nonumber\\
&&\qquad +\int_0^tF(\hat{\eta}_s)\int_{S\times [0,\infty )}\delta 
(x,\eta_s)f_r(x,r)\hat{\eta}_s(dx,dr)ds,\label{ito1}\\
&=&\hspace{-2mm}F(\hat{\eta}_0)+\int_{S\times [0,t]\times [0,\infty )^2}F(\hat{
\eta}_{s-})(e^{-f(x,r)}-1)\one_{[0,\lambda (x,\eta_{s-})]}(u)\tilde {
N}(dx,ds,dr,du)\nonumber\\
&&\qquad +\int_0^tF(\hat{\eta}_s)\Big(\int_{S\times [0,\infty )}\lambda 
(x,\eta_s)(e^{-f(x,r)}-1)e^{-r}\beta (dx)dr\nonumber\\
&&\qquad\qquad\qquad\qquad +\int_{S\times [0,\infty )}\delta (x,\eta_
s)f_r(x,r)\hat{\eta}_s(dx,dr)\Big)ds,\nonumber\end{eqnarray}
where 
$\tilde {N}(dx,ds,dr,du)=N(dx,ds,dr,du)-\beta (dx)\times ds\times 
e^{-r}dr\times du$.  
It follows from \reff{intest} that the stochastic integral 
term on the right is a local martingale, and since $f^{*}$ 
given by (\ref{absderiv}) is in $\hat {{\mathcal C}}$, it follows that  
\begin{equation}\int_0^t\int_{S\times [0,\infty )}\delta (x,\eta_
s)|f_r(x,r)|\hat{\eta}_s(dx,dr)ds<\infty ,\quad t>0.\label{dthint}\end{equation}
Consequently, 
defining 
\begin{eqnarray}
\hat {A}F(\hat{\zeta })&=&F(\hat{\zeta })\Big(\int_{S\times [0,\infty 
)}\lambda (x,\zeta )(e^{-f(x,r)}-1)e^{-r}\beta (dx)dr\nonumber\\
&&\qquad\qquad +\int_{S\times [0,\infty )}\delta (x,\zeta )f_r(x,
r)\hat{\zeta }(dx,dr)\Big),\label{gendef}\end{eqnarray}
any solution of \reff{stoch3a} must be a solution of the 
local martingale problem for $\hat {A}$.  We say that $\hat{\eta}$ is a 
solution of the {\em local martingale problem\/} for $\hat {A}$ if there 
exists a filtration $\{{\mathcal F}_t\}$ such that $\hat{\eta}$ is $\{
{\mathcal F}_t\}$-adapted and 
\begin{equation}M_F(t)=F(\hat{\eta}_t)-F(\hat{\eta}_0)-\int_0^t\hat {
A}F(\hat{\eta}_s)ds\label{lmp}\end{equation}
is a $\{{\mathcal F}_t\}$-local martingale for each $F\in {\mathcal D}(\hat {
A})$, that is, for 
each $F$ of the form $F(\hat{\zeta })=e^{-\int fd\hat{\zeta}}$, $
f\in\hat {{\mathcal C}}$, $f\geq 0$.  In particular, 
let $\bar {f}(x)=\sup_rf(x,r)$ and 
\[\tau_{f,c}=\inf\{t:\int_0^t\int_S\bar {f}(x)\lambda (x,\eta_s)\beta 
(dx)ds>c\}.\]
Then $M_F(\cdot\wedge\tau_{f,c})$ is a martingale.  Note that $\tau_{
f,c}$ is a 
$\{{\mathcal F}^{\eta}_t\}$-stopping time.

Conversely, if $\hat{\eta}$ is a solution of the local martingale 
problem for $\hat {A}$  with sample paths in $D_{\hat {{\mathcal S}}}
[0,\infty )$, then 
under Condition \ref{bound}, (\ref{intest}) and (\ref{dthint}) 
hold.  If $\gamma\in C_{{\Bbb R}}[0,\infty )$ has compact support and  
$f(x,r)=\int_0^r\gamma (u)duc_k(x)$, then $f\in\hat {{\mathcal C}}$ and it follows that 
\[\int_0^t\int_Sc_k(x)\delta (x,\eta_s)\hat{\eta}_s(dx,dr)ds<\infty 
,\quad t>0.\]
To formulate the main theorem of this section, we need 
to introduce the notion of a {\em weak solution\/} of a 
stochastic equation.  

\begin{definition}
{\rm A stochastic process $\tilde{\eta}$ with sample paths in $D_{\hat {
{\mathcal S}}}[0,\infty )$ is 
a {\em weak solution\/} of \reff{stoch2a} if there exists a 
probability space $(\Omega ,{\mathcal F},P)$, a Poisson random measure $
N$ 
on $S\times [0,\infty )^3$ with mean measure $\beta (dx)\times ds
\times e^{-r}dr\times du$ 
and a stochastic process $\hat{\eta}$ defined on $(\Omega ,{\mathcal F}
,P)$, such that 
$\tilde{\eta}$ and $\hat{\eta}$ have the same distribution on $D_{
\hat {{\mathcal S}}}[0,\infty )$, $\hat{\eta}$ is 
adapted to a filtration with respect to which $N$ is 
compatible, and $N$ and $\hat{\eta}$ satisfy \reff{stoch2a}.  }
\end{definition}

\begin{theorem}\label{equivthm}
Suppose that $\lambda$ and $\delta$ satisfy Conditions \ref{bound} and 
\ref{contcnd}.  Then each solution of the stochastic 
equation \reff{stoch2a} (or equivalently, \reff{stoch3a}) 
is a solution of the local martingale problem for $\hat {A}$ 
defined by \reff{gendef}, and
 each solution of the local 
martingale problem for $\hat {A}$ is a weak solution of the 
stochastic equation.  
\end{theorem}

\begin{proof} The first part of the theorem follows from the 
discussion above.  

To prove the second part, we apply a Markov mapping 
result of
Kurtz (1998).  Let $\{D_i\}\subset {\mathcal B}(S\times [0,\infty )^2
)$ be countable, closed 
under intersections, generate ${\mathcal B}(S\times [0,\infty )^2)$, and satisfy 
$\int_{D_i}\beta (dx)e^{-r}\,dr\,ds<\infty$.  Then $N$ is completely determined 
by $N(D_i,t)$.  Define 
\[Z_i(t)=Z_i(0)(-1)^{N(D_i,t)},\]
where $Z_i(0)$ is $\pm 1$.  Note that 
\begin{equation}N(D_i,t)=-\frac 12\int_0^tZ_i(s-)dZ_i(s),\label{cdef}\end{equation}
and if the $Z_i(0)$ are iid with 
$P\{Z_i(0)=1\}=P\{Z_i(0)=-1\}=\frac 12$ and independent of $N$, then 
for each $t\geq 0$, the $Z_i(t)$ are iid and independent of $N$.  
For $z\in \{-1,1\}^{\infty}$, we will let $(-1)^{{\bf 1}_D(x,r,u)}
z$ denote 
\[(-1)^{{\bf 1}_D(x,r,u)}z=((-1)^{{\bf 1}_{D_1}(x,r,u)}z_1,(-1)^{
{\bf 1}_{D_2}(x,r,u)}z_2,\ldots ).\]
Then $Z=(Z_1,Z_2,\ldots )$ is a solution of the martingale 
problems for 
\[CG(z)=\int_{S\times [0,\infty )^2}(G((-1)^{{\bf 1}_D(x,r,u)}z)-
G(z))e^{-r}dr\,du\beta (dx).\]
We can take the domain for $C$ to be the collection of 
functions that depend on only finitely many coordinates 
of $z$.  With this domain, the martingale problem for $C$ is 
well-posed.  

If $\hat{\eta}$ is a solution of \reff{stoch3a}, then $(\hat{\eta }
,Z)$ is a 
solution of the local martingale problem for 
\begin{eqnarray}
\lefteqn{\hat{\mathbb A}(FG)(\hat{\zeta },z) \, } \nonumber \\
&=&F(\hat{\zeta })\Big(\int_{S
\times [0,\infty )^2}\Big(({\bf 1}_{[0,\lambda (x,\zeta )]}(u)e^{
-f(x,r)}+{\bf 1}_{(\lambda (x,\zeta ),\infty )}(u))G((-1)^{{\bf 1}_
D(x,r,u)}z) \nonumber \\
& & \qquad\qquad\qquad\qquad\qquad\qquad\qquad\qquad\qquad\qquad\qquad\,-\,
G(z)\Big)e^{-r}\beta (dx) 
dr\,du\nonumber\\
&&\qquad\qquad\qquad -G(z)\int_{S\times [0,\infty )}\delta (x,\zeta 
)f_r(x,r)\hat{\zeta }(dx,dr)\Big).\label{gendef2}\end{eqnarray}

Let $\tilde{\eta}$ be a solution of the local martingale problem for 
$\hat {A}$.  For $a=(a_1,a_2,\ldots )$ with $a_k>0$, $k=1,2,\ldots$,
define 
\begin{eqnarray*}
\lefteqn{\tau_a(t)\,=}  \\
& & \hspace{-3mm}\inf\{u:\int_0^u1\vee\sum_{k=1}^{\infty}a_k\left[\int_
Sc_k(x)\lambda (x,\tilde\eta_s)\beta (dx)+\int_{S\times S}c_k(x)\delta 
(x,\tilde\eta_s)\tilde\eta_s(dx)\right]ds\geq t\},
\end{eqnarray*}
\[H_a(\zeta )=1\vee\sum_{k=1}^{\infty}a_k\left[\int_Sc_k(x)\lambda 
(x,\zeta )\beta (dx)+\int_{S\times S}c_k(x)\delta (x,\zeta )\zeta 
(dx)\right],\]
and $\tilde{\eta}^a_t=\tilde{\eta}_{\tau_a(t)}$.  Then $\tilde{\eta}^
a$ is a solution of the martingale 
problem for 
\begin{equation}\hat {A}^a\equiv\frac 1{H_a}\hat {A}.\label{truncgen}\end{equation}
For $F\in\!{\mathcal D}(\hat {A})$, $\hat {A}^aF$ is bounded, and we can 
select $a^n=(a_1^n,a_2^n,\ldots )$ so that $a^n\geq a^{n+1}$ 
and $\tau_{a^n}(t)\rightarrow t$ a.s.  

Let $\mu (dz)=\prod_{k=1}^{\infty}(\frac 12\delta_{\{-1\}}(dz_k)+\frac 
12\delta_{\{1\}}(dz_k))$, and set 
$c_G=\int Gd\mu$.  Then 
\[\int\hat {{\Bbb A}}(FG)(\hat{\zeta },z)\mu (dz)=c_G\hat {A}F(\hat{
\zeta }),\]
and more generally, 
\[\int\frac 1{H_a}\hat {{\Bbb A}}A(FG)(\hat{\zeta },z)\mu (dz)=c_
G\frac 1{H_a}\hat {A}F(\hat{\zeta }).\]
Applying Corollary 3.5 of Kurtz (1998) to $H_a^{-1}\hat {{\Bbb A}}$ for 
each $a$, we conclude that if $\tilde{\eta}$ is a solution of the local 
martingale problem for $\hat {A}$, then there exists a solution 
$(\hat{\eta },Z)$ of the local martingale problem for $\hat {{\Bbb A}}$ such 
that $\hat{\eta}$ and $\tilde{\eta}$ have the same distribution.  Finally, 
applying \reff{cdef}, we can construct the corresponding 
Poisson random measure $N$ and show that $\hat{\eta}$ and $N$ 
satisfy \reff{stoch3a}.  \hfill \end{proof}

The natural (local) martingale problem for $\eta$ is really the 
martingale problem for $A$ given by \reff{gener}; 
however, there 
will be solutions $\hat{\eta}$ of the local martingale problem for $
\hat {A}$ 
(and hence of the stochastic equation) such that the 
corresponding $\eta$ is not a solution of the local martingale 
problem for $A$.  Intuitively, conditioned on 
${\mathcal F}_t^{\eta}=\sigma (\eta_s:s\leq t)$, the residual clock times should be 
independent unit exponentials, independent of ${\mathcal F}_t^{\eta}$.  That 
need not be the case, since we are free to pick the 
residual clock times at time zero in any way we please.  
It also need not be the case if the solution of the 
martingale problem fails to be unique.  The following 
results clarify the relationship between the martingale 
problems for $A$ and $\hat {A}$.  

\begin{proposition}\label{projprop} 
Suppose that $\lambda$ and $\delta$ satisfy Conditions \ref{bound} and 
\ref{contcnd}.  
If $\hat{\eta}$ is a solution of the local martingale 
problem for $\hat {A}$ and  at each time $t$, the residual clock 
times are independent of ${\mathcal F}^{\eta}_t$ and are independent unit 
exponentials, then $\eta$ is a solution of the local 
martingale problem for $A$.  
\end{proposition}

\begin{proof}\ By assumption, we can write $\hat{\eta}_t=\sum_i\delta_{
(X_i(t),R_i(t))}$, 
where the $R_i(t)$ are independent unit exponentials, 
independent of ${\mathcal F}^{\eta}_t$, and in particular, independent of $
\eta_t$.  
For $f\in\hat {{\mathcal C}}$ and $F(\hat{\zeta })=e^{-\int_{S\times 
[0,\infty )}f(x,r)\hat{\zeta }(dx,dr)}$, since 
(\ref{lmp}) can be localized by $\{{\mathcal F}_t^{\eta}\}$-stopping times, it 
follows that
\[E[F(\hat{\eta}_t)|{\mathcal F}_t^{\eta}]-E[F(\hat{\eta}_0|{\mathcal F}^{
\eta}_0]-\int_0^tE[\hat {A}F(\eta_s)|{\mathcal F}_s^{\eta}]ds\]
is a $\{{\mathcal F}_t^{\eta}\}$-local martingale.  By the independence of the 
$R_i(t)$, 
\[E[F(\hat{\eta}_t)|{\mathcal F}_t^{\eta}]=\prod_i\int_0^{\infty}e^{-
f(X_i(t),r)}e^{-r}dr=e^{-\int_Sg(x)\eta_t(dx)}\equiv G(\eta_t),\]
where $g$ is defined so that $e^{-g(x)}=\int_0^{\infty}e^{-f(x,r)}
e^{-r}dr$.  
Integrating by parts gives 
\[\int_0^{\infty}e^{-f(x,r)}f_r(x,r)e^{-r}dr=1-\int_0^{\infty}e^{
-f(x,r)}e^{-r}dr=1-e^{-g(x)},\]
and hence 
\begin{eqnarray*}
E[\hat {A}F(\eta_s)|{\mathcal F}_s^{\eta}]&=&G(\eta_s)\int_{S\times [
0,\infty )}\lambda (x,\eta_s)(e^{-f(x,r)}-1)e^{-r}\beta (dx)dr\\
&&\qquad +\sum_j\left(\prod_{i\neq j}\int_0^{\infty}e^{-f(X_i(s),
r)}e^{-r}dr\right)\\
&&\qquad\qquad\int_0^{\infty}e^{-f(X_j(s),r)}\delta (X_j(s),\eta_
s)f_r(X_j(s),r)e^{-r}dr\\
&=&G(\eta_s)\int_S\lambda (x,\eta_s)(e^{-g(x)}-1)\beta (dx)\\
&&\qquad +\sum_j\left(\prod_{i\neq j}e^{-g(X_i(s))}\right)\delta 
(X_j(s),\eta_s)\left(1-e^{-g(X_j(t))}\right)\\
&=&AG(\eta_s),\end{eqnarray*}
and the proposition follows.\hfill \end{proof}

We have the following converse for the previous 
proposition.  

\begin{theorem}\label{projthm} 
Suppose that $\lambda$ and $\delta$ satisfy Conditions \ref{bound} and 
\ref{contcnd}.  
If $\eta$ is a solution of the local martingale problem for $A$, 
then there exists a solution $\hat{\eta}$ of the local martingale 
problem for $\hat {A}$ such that $\eta$ and $\hat{\eta }(\cdot\times 
[0,\infty ))$ have the same 
distribution on $D_{{\mathcal S}}[0,\infty )$ and at each time $t\geq 
0$, the 
residual clock times are independent, unit exponentials 
that are independent of ${\mathcal F}^{\eta}_t$. 
\end{theorem}

\begin{proof}
For $\zeta =\sum_i\delta_{x_i}\in {\mathcal S}$, let $\alpha (\zeta ,
d\hat{\zeta })$ denote the 
distribution on $\hat {{\mathcal S}}$ of $\sum_i\delta_{(x_i,\tau_i)}$, where the $
\tau_i$ are 
independent, unit  exponential random variables.  Then, 
by the calculation in the proof of Proposition 
\ref{projprop}, 
\[G(\zeta )=\int_{\hat {{\mathcal S}}}F(\hat{\zeta })\alpha (\zeta ,d
\hat{\zeta })\qquad AG(\zeta )=\int_{\hat {{\mathcal S}}}\hat {A}F(\hat{
\zeta })\alpha (\zeta ,d\hat{\zeta }),\]
 for $F\in {\mathcal D}(\hat {A})$.  More generally, 
$A^aG(\zeta )=\int_{\hat {{\mathcal S}}}\hat {A}^aF(\hat{\zeta })\alpha 
(\zeta ,d\hat{\zeta })$, where $\hat {A}^a$ is defined as in 
\reff{truncgen}.  The theorem then follows by Corollary 
3.5 of Kurtz (1998).\hfill 
\end{proof}

\begin{corollary}
Let $\nu\in {\mathcal P}({\mathcal S})$, and define $\hat{\nu}\in {\mathcal P}
(\hat {{\mathcal S}})$ by 
\[\int_{\hat {{\mathcal S}}}hd\hat{\nu }=\int_{{\mathcal S}}\int_{\hat {{\mathcal S}}}
h(\hat{\zeta })\alpha (\zeta ,d\hat{\zeta })\nu (d\zeta ).\]
  If uniqueness holds for 
the martingale problem for $(\hat {A},\hat{\nu })$, or equivalently, weak 
uniqueness holds for the stochastic equation 
\reff{stoch3a}, then uniqueness holds for the martingale 
problem for $(A,\nu )$.  
\end{corollary}

\begin{proof}
If $\eta$ is a solution of the martingale problem for 
$(A,\nu )$, then Theorem \ref{projthm} gives a corresponding 
solution of the martingale problem for $(\hat {A},\hat{\nu })$.
Uniqueness for the latter then implies uniqueness of the 
former.\hfill \end{proof}

\subsection{Existence}
We now turn to the question of existence of solutions of 
(\ref{stoch1a}).  We assume that Conditions \ref{bound}\  
and \ref{contcnd} hold.  The pair $(\lambda ,\delta )$ will be called 
{\em attractive\/} if $\zeta_1\subset\zeta_2$ implies $\lambda (x
,\zeta_1)\leq\lambda (x,\zeta_2)$ and 
$\delta (x,\zeta_1)\geq\delta (x,\zeta_2)$.  If $(\lambda ,\delta 
)$ is attractive and we set 
$\eta^0\equiv 0$, then $\eta^n$ defined by
\begin{eqnarray}
\eta^{n+1}_t(B)&=&\int_{B\times [0,t]\times [0,\infty )^2}\one_{[
0,\lambda (x,\eta^n_{s-})]}(u)\one_{(\int_s^t\delta (x,\eta^n_v)\,
dv,\infty )}(r)N(dx,ds,dr,du)\nonumber\\
&&\quad\quad\quad\quad\quad\quad\quad\quad +\int_{B\times [0,\infty 
)}\one_{(\int_0^t\delta (x,\eta^n_s)\,ds,\infty )}(r)\hat{\eta}_0
(dx,dr)\label{stocheqe}\end{eqnarray}
is monotone increasing and either $\eta^n$ converges to a 
process with values in ${\mathcal S}$, or 
\begin{equation}\int_0^T\int_Sc_k(x)\lambda (x,\eta^n_s)\beta (dx
)ds\rightarrow\infty ,\label{diverg}\end{equation}
for some $T$ and $k$.  To see this, let
\[\tau_c^n=\inf\{t:\int_0^t\int_Sc_k(x)\lambda (x,\eta^n_s)\beta 
(dx)ds>c\}.\]
Then
\begin{eqnarray*}
&&E[\sup_{t\leq T\wedge\tau_c^n}\left(\int_Sc_k(x)\eta_s^{n+1}(dx
)-\int_{A\times [0,\infty )}\one_{(\int_0^t\delta (x,\eta^n_s)\,d
s,\infty )}(r)\hat\eta_0(dx,dr)\right)]\\
&&\qquad\leq E[\int_0^{T\wedge\tau_c^n}\int_Sc_k(x)\lambda (x,\eta^
n_s)\beta (dx)ds]\\
&&\qquad\leq c,\end{eqnarray*}
and $\tau_c^1\geq\tau_c^2\geq\cdots$.  Either 
\begin{equation}\lim_{c\rightarrow\infty}\lim_{n\rightarrow\infty}
\tau_c^n=\infty\label{tauinf}\end{equation}
or 
(\ref{diverg}) holds for some $T$.

If (\ref{tauinf}) holds almost surely, the limit $\eta^{\infty}$ is the 
minimal solution of (\ref{stoch1a}) in the sense that any 
other solution $\eta$ will satisfy $\eta^{\infty}_t(B)\leq\eta_t(
B)$ for all 
$B\in {\mathcal B}(S)$ and $t\geq 0$.  

For an arbitrary pair $(\lambda ,\delta )$ satisfying Conditions \ref{bound} and 
\ref{contcnd}, we define an attractive pair by setting
\[\bar{\lambda }(x,\zeta )=\sup_{\zeta'\subset\zeta}\lambda (x,\zeta'
)\qquad\underline {\delta}(x,\zeta )=\inf_{\zeta'\subset\zeta}\delta 
(x,\zeta').\]
Let $\eta_0$ be an 
${\mathcal S}$-valued random variable independent of $N$, and let $\hat{
\eta}_0$ be 
defined as in (\ref{etahat}). 
We assume that $\bar{\lambda}$ satisfies (\ref{brest}), which implies
\begin{equation}\int c_k(x)\bar{\lambda }(x,\zeta )\beta (dx)<\infty 
,\quad\zeta\in S,k=1,2,\ldots ,\label{atbnd}\end{equation}
 and that
there exists a solution $\bar{\eta}$ for the 
pair $(\bar{\lambda },\underline {\delta})$.   

We consider a different sequence of approximate 
equations.  Let $\{K_n\}$ be the sets in the definition of ${\mathcal C}$, 
and let $\eta^n$ satisfy
\begin{eqnarray}
\lefteqn{\eta^n_t(B) \, =} \nonumber \\
&&\hspace{-2mm}\int_{B\times [0,t]\times [0,\infty )^2}\one_{[0,\lambda 
(x,\bar{\eta}_{s-}\cap K_n\cap\eta^n_{s-})]}(u)\one_{(\int_s^t\delta 
(x,\bar{\eta}_v\cap K_n\cap\eta^n_v)\,dv,\infty )}(r)N(dx,ds,dr,d
u)\nonumber\\
&&\quad\quad\quad\quad\quad\quad\quad\quad +\int_{B\times [0,\infty 
)}\one_{(\int_0^t\delta (x,\bar{\eta}_v\cap K_n\cap\eta^n_v)\,dv,
\infty )}(r)\hat{\eta}_0(dx,dr).\label{stocheqe2}\end{eqnarray}
Existence and uniqueness for this equation follow from 
the fact that only finitely many births can occur in a 
bounded time interval in $K_n$.  Consequently, the equation 
can be solved from one such birth to the next.
Since $\lambda (x,\bar{\eta}_{s-}\cap K_n\cap\eta^n_{s-})\leq\bar{
\lambda }(x,\bar{\eta}_{s-})$  and 
$\delta (x,\bar{\eta}_s\cap K_n\cap\eta^n_v)\geq\underline {\delta}
(x,\bar{\eta}_s)$, it follows that $\eta^n_t\subset\bar{\eta}_t$ and hence 
that 
\begin{eqnarray}
\eta^n_t(B)&=&\int_{B\times [0,t]\times [0,\infty )^2}\one_{[0,\lambda 
(x,K_n\cap\eta^n_{s-})]}(u)\one_{(\int_s^t\delta (x,K_n\cap\eta^n_
v)\,dv,\infty )}(r)N(dx,ds,dr,du)\nonumber\\
&&\quad\quad\quad\quad\quad\quad\quad\quad +\int_{B\times [0,\infty 
)}\one_{(\int_0^t\delta (x,K_n\cap\eta^n_v)\,dv,\infty )}(r)\hat{
\eta}_0(dx,dr).\label{stocheqe2b}\end{eqnarray}
Also, note that for $g\in {\mathcal C}$,
\begin{equation}\int_0^t\int_Sg(x)\delta (x,K_n\cap\eta^n_s)\eta^
n_s(dx)ds\leq\int_0^tg(x)r\one_{[0,\bar{\lambda }(x,\bar{\eta}_{s
-})]}(u)N(dx,ds,dr,du)<\infty .\label{dthest}\end{equation}
Define 
 $F(\hat{\zeta })=e^{-\int_{S\times [0,\infty )}f(x,r)\hat{\zeta }
(dx,dr)}$, $f\in\hat {{\mathcal C}}$ nonnegative.  Setting
\begin{eqnarray}
\hat {A}_nF(\hat{\zeta })&=&F(\hat{\zeta })\Big(\int_{S\times [0,
\infty )}\lambda (x,K_n\cap\zeta )(e^{-f(x,r)}-1)e^{-r}\beta (dx)
dr\nonumber\\
&&\qquad\qquad +\int_{S\times [0,\infty )}\delta (x,K_n\cap\zeta 
)f_r(x,r)\hat{\zeta }(dx,dr)\Big),\label{apgendef}\end{eqnarray}
as in (\ref{ito1}),
\[F(\hat{\eta}^n_t)-F(\hat{\eta}_0)-\int_0^t\hat {A}_nF(\hat{\eta}^
n_s)ds\]
is a local martingale.  

Uniqueness for (\ref{stocheqe2b}) implies that the 
residual clock times at time $t$ are conditionally independent, 
unit exponentials given ${\mathcal F}^{\eta}_t$.  Consequently, as in 
Proposition \ref{projprop}, for $G(\zeta )=e^{-\int_Sg(x)\zeta (d
x)}$, $g\in {\mathcal C}$ 
nonnegative, 
and
\[A_nG(\zeta )=\int (G(\zeta +\delta_x)-G(\zeta ))\lambda (x,K_n\cap
\zeta )\beta (dx)+\int (G(\zeta -\delta_x)-G(\zeta ))\delta (x,K_
n\cap\zeta )\zeta (dx),\]
\begin{equation}G(\eta^n_t)-G(\eta_0^t)-\int_0^tA_nG(\eta^n_s)ds\label{lmgex}\end{equation}
is a local martingale.
Exploiting the fact that $\eta^n_t\subset\bar{\eta}_t$, 
the relative compactness of $\{\eta^n\}$, in the sense of 
convergence in distribution in $D_{{\mathcal C}}[0,\infty )$ follows.  

\begin{proposition}Suppose that 
 Conditions \ref{bound}  
and \ref{contcnd} hold.
If 
$(x,\zeta )\rightarrow\lambda (x,\zeta )$ and $(x,\zeta )\rightarrow
\delta (x,\zeta )$ are continuous on $S\times {\mathcal C}$, 
then $\zeta\rightarrow AG(\zeta )$ is  continuous, and  
any limit point of $\{\eta^n\}$ is a 
solution of the local martingale problem for $A$, and hence 
a weak solution of  \reff{stoch2a}.  
\end{proposition}
 
\begin{proof}
By (\ref{atbnd}), we can select $a_k$ so that 
\[\Gamma (t)\equiv\int_0^t1\vee\sum_ka_k\int_Sc_k(x)\bar{\lambda }
(x,\bar{\eta}_s)\beta (dx)ds<\infty ,\quad\forall t>0\quad a.s.\]
and by (\ref{dthest}), it follows that 
\[\tau_m=\inf\{t:\Gamma (t)\geq m\}\]
is a localizing sequence for (\ref{lmgex}) for all $g$ and $n$.  
The estimates also give the necessary uniform 
integrability to ensure that limit points of (\ref{lmgex}) 
are local martingales.
\end{proof}

\subsection{Existence and Uniqueness}

If $\sup_{\zeta\in {\mathcal S}}\int_S\lambda (x,\zeta )\beta
(dx)<\infty$, then a solution of  
\reff{stoch1a} has only finitely many births per unit 
time and it is easy to see that \reff{stoch1a} has a 
unique solution.  Condition \ref{bound}, however, only 
ensures that there are finitely many births per unit 
time in each $K_k$, and uniqueness requires 
additional conditions.  The conditions we use are 
essentially the same as those used for existence and 
uniqueness of the solution of the time change system in 
Garcia (1995).  From now on, we are going to assume 
that $\delta (x,\eta )=1$, for all $x\in S$ and $\eta\in {\mathcal S}$.  

Let $N$ be a Poisson random measure on $S\times [0,\infty )^3$ with 
mean measure $\beta (dx)\times ds\times e^{-r}dr\times du$.  Let $
\eta_0$ be an 
${\mathcal S}$-valued random variable independent of $N$, and let $\hat{
\eta}_0$ be 
defined as in (\ref{etahat}).  Suppose $\{{\mathcal F}_t\}$ is a filtration such that 
$\hat{\eta}_0$ is ${\mathcal F}_0$-measurable and $N$ is $\{{\mathcal F}_
t\}$-compatible.  We 
consider the equation 
\begin{eqnarray}
\eta_t(B)&=&\int_{B\times [0,t]\times [0,\infty )^2}\one_{[0,\lambda 
(x,\eta_{s-})]}(u)\one_{(t-s,\infty )}(r)N(dx,ds,dr,du)\nonumber\\
&&\quad\quad\quad\quad\quad\quad\quad\quad +\int_{B\times [0,\infty 
)}\one_{(t,\infty )}(r)\hat{\eta}_0(dx,dr).\label{stoch5}\end{eqnarray}

\begin{theorem}\ 
\label{thmeu}Assume Conditions \ref{bound} and 
\ref{contcnd}.  Suppose that 
$$a(x,y)\geq\sup_{\eta}|\lambda (x,\eta +\delta_y)-\lambda (x,\eta 
)|$$ and that there exists a 
positive function $c$ such that 
\[M=\sup_x\int_S\frac {c(x)a(x,y)}{c(y)}\beta (dy)<\infty .\]
Then, there exists a unique solution of \reff{stoch5}.  
\end{theorem}\ 

\begin{example}{\rm Let $d(x,\eta )=\inf\{d_S(x,y):y\in\eta \}$, where
    $d_S$ is a distance in $S$ such that $(S,d_S)$ is complete
    separable metric space.  Suppose $\lambda (x,\eta )=h(d(x,\eta
    ))$. Then $a(x,y)=\sup_{r>d_S(x,y)}|h(r)-h(d_S(x,y))|$.  If $h$ is
    increasing, then $a(x,y)=h(\infty )-h(d_S(x,y))$ and
\[|\lambda (x,\eta^1)-\lambda (x,\eta^2)|\leq\int (h(\infty )-h(d_
S(x,y)))|\eta^1-\eta^2|(dy).\]
  If $h$ is decreasing, then $a(x,y)=h(d_S(x,y))-h(\infty )$ and
\[|\lambda (x,\eta^1)-\lambda (x,\eta^2)|\leq\int (h(d_S(x,y))-h(
\infty ))|\eta^1-\eta^2|(dy).\]
 }
\end{example}

Theorem \ref{thmeu}\ is a consequence of the following lemmas 
that hold under the conditions of the theorem.

\begin{lemma}\ 
\label{lemmaTV} For any $\eta^1,\eta^2\in {\mathcal S}$  we have 
\begin{equation}|\lambda (x,\eta^1)-\lambda (x,\eta^2)|\le\int_Sa
(x,y)\,|\eta^1-\eta^2|(dy).\label{eqcolc1}\end{equation}
\end{lemma}

\begin{proof}
Since $\eta^1$ and $\eta^2$ contain countably many points, 
there exist $\{y_1,y_2,\ldots \}$ and $\{z_1,z_2,\ldots \}$ such that 
\[\eta^2=\eta^1+\sum_{i=1}^I\delta_{y_i}-\sum_{j=1}^J\delta_{z_j}\]
(where $I$ and $J$ may be infinity) and hence 
\[|\eta^1-\eta^2|(B)=\sum_{i=1}^I\delta_{y_i}(B)+\sum_{j=1}^J\delta_{
z_j}(B).\]
By the definition of $a$ and Condition \ref{contcnd} 
\begin{eqnarray}
|\lambda (x,\eta^1)-\lambda (x,\eta^2)|&=&\lim_{n\rightarrow\infty}
|\lambda (x,\eta^1)-\lambda (x,\eta^1+\sum_{i=1}^{I\wedge n}\delta_{
y_i}-\sum_{j=1}^{J\wedge n}\delta_{z_j})|\label{lamest}\\
&\le&\lim_{n\rightarrow\infty}\left(\sum_{i=1}^{I\wedge n}a(x,y_i
)+\sum_{j=1}^{J\wedge n}a(x,z_j)\right)\nonumber\\
&\le&\int_Sa(x,y)\,|\eta^1-\eta^2|(dy).\nonumber\end{eqnarray}
 
\hfill \end{proof}

Define 
\[\eta_0(B,t)=\int_{B\times [0,\infty )}\one_{(t,\infty )}(r)\hat{
\eta}_0(dx,dr).\]
Let $\eta$ be $\{{\mathcal F}_t\}$-adapted with sample paths in $D_{{\mathcal S}}
[0,\infty )$.  
Then by Condition \ref{bound} 
\begin{equation}\label{c302}\Phi\eta_t(B)=\eta_0(B,t)+\int_{B\times 
[0,t]\times [0,\infty )^2}\one_{[0,\lambda (x,\eta_s)]}(u)\one_{(
t-s,\infty )}(r)N(dx,ds,dr,du)\end{equation}
defines a process adapted to $\{{\mathcal F}_t\}$ with sample paths in 
$D_{{\mathcal S}}[0,\infty )$.

\begin{lemma}\label{lemma29} 
Let $\eta^1$ and $\eta^2$ be adapted to  $\{{\mathcal F}_t\}$ and have sample paths 
in $D_{{\mathcal S}}[0,\infty$).  Then 
\begin{eqnarray}
&&\sup_xc(x)E[\int_Sa(x,y)|\Phi\eta^1(t)-\Phi\eta^2(t)|(dy)]
\label{eq7}\\
&&\qquad\le M\int_0^t\sup_x\,c(x)E[\int_Sa(x,y)|\eta^1_s-\eta^2_s
|(dy)]e^{-(t-s)}\,ds.\nonumber\end{eqnarray}
\end{lemma}

\begin{proof}
Let $\xi^i=\Phi\eta^i$.  Then 
\begin{eqnarray}
&&\sup_zc(z)E[\int_Sa(z,x)|\xi^1_t-\xi^2_t|(dx)]\nonumber\\
&&\quad\le\sup_zc(z)E[\int_{S\times [0,t]\times [0,\infty )^2}a(z
,x)|\one_{[0,\lambda (x,\eta^1_s)]}(u)-\one_{[0,\lambda (x,\eta^2_
s)]}(u)|\nonumber\\
&&\qquad\qquad\qquad\qquad\qquad\qquad\qquad\qquad\qquad\one_{(t-
s,\infty )}(r)N(dx,ds,dr,du)]\nonumber\\
&&\quad\le\sup_zc(z)E[\int_{S\times [0,t]}a(z,x)|\lambda (x,\eta^
1_s)-\lambda (x,\eta^2_s)|e^{-(t-s)}\beta (dx)ds]\nonumber\\
&&\quad\le\sup_zc(z)\,\int_{S\times [0,t]}a(z,x)E[\int_Sa(x,y)|\eta_
s^1-\eta^2_s|(dy)]\,e^{-(t-s)}\,\beta (dx)ds\nonumber\\
&&\quad\le\sup_zc(z)\,\int_S\frac {a(z,x)}{c(x)}\beta (dx)\,\int_
0^t\sup_x\,c(x)E[\int_Sa(x,y)|\eta^1_s-\eta^2_s|(dy)]e^{-(t-s)}\,
ds\nonumber\\
&&\quad\le M\,\int_0^t\sup_x\,c(x)E[\int_Sa(x,y)|\eta^1_s-\eta^2_
s|(dy)]e^{-(t-s)}\,ds.\label{lipest}\end{eqnarray}
\hfill \end{proof}

\begin{proof}\ (Theorem \ref{thmeu}) Uniqueness follows by 
\reff{eq7} and Gronwall's inequality.  To prove existence, 
we proceed by iteration.  Let $\eta^0_t=\eta_0(\cdot ,t)$, and for  $
n\ge 1$, 
define $\eta^{n+1}=\Phi\eta^n$.  Then 
\begin{eqnarray}
&&\sup_xc(x)E[\int_Sa(x,y)|\eta^{n+1}_t-\eta^n_t|(dy)]\nonumber\\
&&\quad\le M\int_0^t\sup_xc(x)E[\int_Sa(x,y)|\eta^n_s-\eta^{n-1}_
s|(dy)]e^{-(t-s)}\,ds\nonumber\\
&&\quad\le M^2\int_0^t\int_0^{s_1}\hspace{-2mm}\sup_xc(x)E[\int_Sa(x,y)|\eta^{
n-1}_{s_2}-\eta_{s_2}^{n-2}|(dy)]e^{-(s_1-s_2)}\,ds_2\,e^{-(t-s_1
)}\,ds_1\nonumber\\
&&\quad\le M^n\int_0^t\int_0^{s_1}\ldots\int_0^{s_n-1}\sup_xc(x)E
[\int_Sa(x,y)|\eta^1_{s_{n-1}}-\eta_{s_{n-1}}^0|(dy)\nonumber\\
&&\qquad\qquad\qquad\qquad\qquad\qquad e^{-(s_{n-1}-s_n)}ds_n\,\ldots\,
e^{-(t-s_1)}\,ds_1.\nonumber\end{eqnarray}
Therefore, there exists $C>0$ such that 
\[\sup_xc(x)E[\int_Sa(x,y)|\eta^{n+1}_t-\eta^n_t|(dy)]\leq\frac {
C^nt^n}{n!}\sup_{s\leq t}\sup_xc(x)E[\int_Sa(x,y)|\eta_s^1-\eta^0_
s|(dy)],\]
and the convergence of $\eta^n$ to a solution of \reff{stoch5} 
follows.  \hfill \end{proof}

\setcounter{equation}{0}

\section{Ergodicity for spatial birth and death 
processes}  \label{sectiondeath} 

\subsection{Temporal ergodicity}
The statement that a Markov process is ergodic can 
carry several meanings.  At a minimum, it means that 
there exists an unique stationary distribution for the 
process. Under this condition, the corresponding 
stationary process is ergodic in the sense of triviality of 
its tail  $\sigma$-algebra.  A second, stronger meaning of 
ergodicity for Markov processes is that for all initial 
distributions, the distribution 
of the process at time $t$ converges to the (unique) stationary 
distribution as $t\rightarrow\infty$.

One approach to the first kind of ergodicity involves 
using the stochastic equation to construct a ``coupling 
form the past.''  Following an idea of Kendall and M\o ller 
(2000),
for 
$\eta^1\subset\eta^2$,
define
\[\bar{\lambda }(x,\eta^1,\eta^2)=\sup_{\eta^1\subset\eta\subset\eta^
2}\lambda (x,\eta )\qquad\underline {\lambda}(x,\eta^1,\eta^2)=\inf_{
\eta^1\subset\eta\subset\eta^2}\lambda (x,\eta ).\]
Note that for $\eta^1\subset\eta^2$ 
\[|\bar{\lambda }(x,\eta^1,\eta^2)-\underline {\lambda}(x,\eta^1,
\eta^2)|\leq\int_Sa(x,y)|\eta^1-\eta^2|(dy).\]

We assume that $N$ is defined on $S\times (-\infty ,\infty )\times 
[0,\infty )^2$, that 
is, for all positive and negative time, and consider a 
system starting from time $-T$, that is, for $t\geq -T$
\begin{eqnarray}
\eta_t^{1,T}(B)&=&\int_{B\times [-T,t]\times [0,\infty )^2}\one_{
[0,\underline {\lambda}(x,\eta^{1,T}_{s-},\eta^{2,T}_{s-}))}(u)\one_{
(t-s,\infty )}(r)N(dx,ds,dr,du)\nonumber\\
&&\quad\quad\quad\quad\quad\quad\quad\quad +\int_{B\times [0,\infty 
)}\one_{(t+T,\infty )}(r)\hat{\eta}_{-T}^{1,T}(dx,dr)\nonumber\\
\eta_t^{2,T}(B)&=&\int_{B\times [0,t]\times [0,\infty )^2}\one_{[
0,\bar{\lambda }(x,\eta^{1,T}_{s-},\eta_{s-}^{2,T})]}(u)\one_{(t-
s,\infty )}(r)N(dx,ds,dr,du)\nonumber\\
&&\quad\quad\quad\quad\quad\quad\quad\quad +\int_{B\times [0,\infty 
)}\one_{(t+T,\infty )}(r)\hat{\eta}_{-T}^{2,T}(dx,dr),\label{coupsys}\end{eqnarray}
where we require $\eta^{1,T}_{-T}\subset\eta^{2,T}_{-T}$.  Suppose $
\lambda (x,\eta )\leq\Lambda (x)$ for all $\eta$ 
and 
\begin{equation}\int_Sc_k(x)\Lambda (x)\beta (dx)<\infty ,\qquad 
k=1,2,\ldots ,\label{biglam}\end{equation}
which implies 
Condition \ref{bound}, and suppose Condition \ref{contcnd} 
holds.
Then we can obtain a 
solution of (\ref{coupsys}) by iterating 
\begin{eqnarray}
\eta_t^{1,T,n+1}(B)&=&\int_{B\times [-T,t]\times [0,\infty )^2}\one_{
[0,\underline {\lambda}(x,\eta^{1,T,n}_{s-},\eta^{2,T,n}_{s-}))}(
u)\one_{(t-s,\infty )}(r)N(dx,ds,dr,du)\nonumber\\
&&\quad\quad\quad\quad\quad\quad\quad\quad +\int_{B\times [0,\infty 
)}\one_{(t+T,\infty )}(r)\hat{\eta}_{-T}^{1,T}(dx,dr)\nonumber\\
\eta_t^{2,T,n+1}(B)&=&\int_{B\times [-T,t]\times [0,\infty )^2}\one_{
[0,\bar{\lambda }(x,\eta^{1,T,n}_{s-},\eta_{s-}^{2,T,n})]}(u)\one_{
(t-s,\infty )}(r)N(dx,ds,dr,du)\nonumber\\
&&\quad\quad\quad\quad\quad\quad\quad\quad +\int_{B\times [0,\infty 
)}\one_{(t+T,\infty )}(r)\hat{\eta}_{-T}^{2,T}(dx,dr),\label{coupsys2}\end{eqnarray}
where we take $\eta^{1,T,1}_t\equiv\emptyset$ and
\begin{eqnarray*}
\eta_t^{2,T,1}(B)&=&\int_{B\times [-T,t]\times [0,\infty )^2}\one_{
[0,\Lambda (x)]}(u)\one_{(t-s,\infty )}(r)N(dx,ds,dr,du)\\
&&\quad\quad\quad\quad\quad\quad\quad\quad +\int_{B\times [0,\infty 
)}\one_{(t+T,\infty )}(r)\hat{\eta}_{-T}^{2,T}(dx,dr).\end{eqnarray*}
Note that $\eta^{1,T,n}\subset\eta^{2,T,n}$, $\{\eta^{1,T,n}\}$ is monotone increasing, 
and $\{\eta^{2,T,n}\}$ is monotone decreasing, and the limit, which 
must exist, will be a solution of (\ref{coupsys}).

For $C\subset\R$, define $(C+t)=\{(s+t):s\in C\}$, and define the 
time-shift of $N$ by
$R_tN(B\times C\times D\times E)=N(B\times (C+t)\times D\times E)$.
Taking $T=\infty$ in (\ref{coupsys2}),  
the iterates
\begin{eqnarray}
\eta_t^{1,\infty ,n+1}(B)\hspace{-2mm}&=&\hspace{-2mm}\int_{B\times (-\infty ,t]\times [0,\infty 
)^2}\hspace{-2mm}\one_{[0,\underline {\lambda}(x,\eta^{1,\infty ,n}_{s-},\eta^{
2,\infty ,n}_{s-})}(u)\one_{(t-s,\infty )}(r)N(dx,ds,dr,du)\nonumber\\
\label{coupsys3}\\
\eta_t^{2,\infty ,n+1}(B)\hspace{-2mm}&=&\hspace{-2mm}\int_{B\times (-\infty ,t]\times [0,\infty 
)^2}\hspace{-2mm}\one_{[0,\bar{\lambda }(x,\eta^{1,\infty ,n}_{s-},\eta_{s-}^{
2,\infty ,n})]}(u)\one_{(t-s,\infty )}(r)N(dx,ds,dr,du),\nonumber\end{eqnarray}
satisfy $\eta^{m,\infty ,n}_t=H^{m,n}(R_tN)$, $m=1,2$, for deterministic 
mappings 
$H^{m,n}:{\mathcal N}(S\times (-\infty ,\infty )\times [0,\infty )^2)
\rightarrow {\mathcal N}(S)$  and the limits $\eta_t^{m,\infty}$ 
satisfy 
\begin{equation}\eta^{m,\infty}_t=H^m(R_tN),\label{map2}\end{equation}
where $H^m=\lim_{n\rightarrow\infty}H^{m,n}$.  It 
follows that $(\eta_t^{1,\infty},\eta_t^{2,\infty})$ is stationary and ergodic.

\begin{lemma}
Suppose that $\lambda$ satisfies (\ref{biglam}) and Condition 
\ref{contcnd}. 
 Then 
\[\eta^{1,\infty}_t\equiv\lim_{n\rightarrow\infty}\eta_t^{1,\infty 
,n}\mbox{\rm \ and }\eta^{2,\infty}_t\equiv\lim_{n\rightarrow\infty}
\eta_t^{2,\infty ,n}\]
 exist and are stationary.
\end{lemma}

Applying 
Theorem \ref{equivthm}, any 
stationary solution of the martingale problem can be 
represented as a weak solution $\eta$ of the stochastic 
equation on the doubly infinite time interval and hence 
coupled to versions of $\eta^{1,\infty ,n}$ and $\eta^{2,\infty ,
n}$ so that 
$\eta^{1,\infty ,n}_t\subset\eta_t\subset\eta^{2,\infty ,n}_t$, $
-\infty <t<\infty$.  Consequently, we have 
the following.

\begin{lemma}\label{uniquesd}
Suppose that $\lambda$ satisfies (\ref{biglam}) and Condition 
\ref{contcnd}.  If 
\[\lim_{n\rightarrow\infty}\int_Sc_k(x)|\eta_t^{2,\infty ,n}-\eta^{
1,\infty ,n}_t|(dx)=0\quad a.s.\]
for $k=1,2,\ldots$, then $\eta\equiv\eta^{2,\infty}=\eta^{1,\infty}$ a.s. is a stationary 
solution of (\ref{stoch5}) and the distribution 
of $\eta_t^{2,\infty}$ is the unique stationary distribution for $
A$.
\end{lemma}

\begin{theorem}\ 
\label{contre} Let $\lambda :S\times {\mathcal N}(S)\rightarrow [0,\infty 
)$ satisfy the conditions of 
Theorem \ref{thmeu} with $M<1$.  
Then $\eta\equiv\eta^{2,\infty}=\eta^{1,\infty}$ a.s. is a stationary 
solution of (\ref{stoch5}) and the distribution 
of $\eta_t^{2,\infty}$ is the unique stationary distribution for $
A$.
\end{theorem}

\begin{proof}
As in the proof of Theorem \ref{thmeu},
\begin{eqnarray}
&&\sup_xc(x)E[\int_Sa(x,y)|\eta^{2,\infty ,n+1}_t-\eta^{1,\infty 
,n+1}_t|(dy)]\nonumber\\
&&\quad\le M\int_{-\infty}^t\sup_xc(x)E[\int_Sa(x,y)|\eta^{2,\infty 
,n}_s-\eta^{1,\infty ,n}_s|(dy)]e^{-(t-s)}\,ds\nonumber\\
&&\quad =M\sup_xc(x)E[\int_Sa(x,y)|\eta^{2,\infty ,n}_t-\eta^{1,\infty 
,n}_t|(dy)],\nonumber\end{eqnarray}
where the equality follows by the stationarity of 
$\eta^{2,\infty ,n}$ and $\eta^{1,\infty ,n}$.  Since the expression on the left is 
nonincreasing, its limit $\rho$ exists, and we have $0\leq\rho\leq 
M\rho$.  
But $M<1$, so $\rho =0$.
\end{proof}

\medskip

\begin{definition} {\rm
$\lambda (x,\cdot )$ is {\em nondecreasing}, if $\eta_1\subset\eta_
2$ implies $\lambda (x,\eta_1)\le\lambda (x,n_2)$.  }
\end{definition}

Note that if $\lambda$ is nondecreasing, then for $\eta_1\subset\eta_
2$, 
$\bar{\lambda }(x,\eta_1,\eta_2)=\lambda (x,\eta_2)$  and $\underline {
\lambda}(x,\eta_1,\eta_2)=\lambda (x,\eta_1)$.  The following 
lemma is immediate.  

\begin{lemma}
Let $\lambda$ be nondecreasing and satisfy (\ref{biglam})  and 
Condition \ref{contcnd}.  Then $\eta^{1,\infty}_t\equiv\lim_{n\rightarrow
\infty}\eta_t^{1,\infty ,n}$ and 
$\eta^{2,\infty}_t\equiv\lim_{n\rightarrow\infty}\eta_t^{2,\infty 
,n}$ are, respectively, the minimal and 
maximal stationary solutions of the martingale problem 
for $A$.
\end{lemma}

For $\lambda$ nondecreasing, the minimal stationary distribution 
can also easily be obtained as a temporal limit.

\begin{lemma}\ \label{mono}
If uniqueness holds for (\ref{stoch5}) and $\lambda (x,\cdot )$ is nondecreasing, 
 then the process $\eta_t$ is attractive, that is 
\begin{equation}\label{r2}\eta_0^1\subset\eta_0^2\quad\mbox{\rm implies}
\quad\eta_t^1\subset\eta_t^2\end{equation}
for all $t\ge 0$.  
\end{lemma}
\begin{proof}
The conclusion is immediate from coupling the two 
processes using the same underlying Poisson random 
measure.  \end{proof}

\begin{theorem}\ Suppose $\lambda$ satisfies (\ref{biglam}) and 
Condition \ref{contcnd}.
If $\lambda (x,\cdot )$ is nondecreasing and $\eta_0=\emptyset$, then the 
distribution of $\eta_t$ converges to the  
minimal stationary 
distribution.
\end{theorem}\ 

\begin{proof}
Note that, if we set $\eta^t_{-t}=\emptyset$, then $\eta_t$ has the same 
distribution as $\eta^t_0$, and by Lemma \ref{mono}, $\eta^t_s\subset
\eta^{1,\infty}_s$ 
for $s\geq -t$.  Since for each $s\geq -t$, $\eta_s^t$ is monotone 
increasing in $t$, $\tilde{\eta}_s=\lim_{t\rightarrow\infty}\eta^
t_s$ exists and must be a 
stationary process.  Since $\eta^{1,\infty}$ is the minimal stationary 
process, we must have $\tilde{\eta}_s=\eta^{1,\infty}_s$.\hfill \end{proof}

The same argument gives the following result on the 
maximal stationary distribution.

\begin{theorem}\ Suppose $\lambda$ satisfies (\ref{biglam}) and 
Condition \ref{contcnd}.
If $\lambda (x,\cdot )$ is nondecreasing and 
\[\eta_0(B)=\int_{B\times (-\infty ,0]\times [0,\infty )^2}\one_{
[0,\Lambda (x)]}(u)\one_{(t-s,\infty )}(r)N(dx,ds,dr,du),\]
then the 
distribution of $\eta_t$ converges to the  
maximal stationary 
distribution.
\end{theorem}\ 

\begin{remark}
{\rm\ Note that $\eta_0$ is a Poisson random measure with 
mean measure $\mu (B)=\int_B\Lambda (x)\beta (dx)$.  }
\end{remark}

\medskip

We can also use the stochastic equation and estimates 
similar to those used in the proof of uniqueness to give 
conditions for ergodicity in the sense of convergence as 
$t\rightarrow\infty$ for all initial distributions.

\begin{theorem}\ 
\label{thte} Let $\lambda :S\times {\mathcal N}(S)\rightarrow [0,\infty 
)$ satisfy the conditions of 
Theorem \ref{thmeu} with $M<1$.  Then the process 
obtained as a solution of the system of  stochastic 
equations \reff{stoch5} is temporally ergodic and the  
rate of convergence is exponential.  
\end{theorem}

\begin{proof}
Suppose $\eta^1$ and $\eta^2$ are solutions of the system 
\reff{stoch5} with distinct initial configurations $\eta^1_0$ and $
\eta^2_0$ 
(equivalently, $\hat{\eta}^1_0$ and $\hat{\eta}^2_0$).  Then, by exactly the same 
argument as used for the proof of Lemma \ref{lemma29} 
we obtain 
\begin{eqnarray}
&&\sup_xc(x)\,E[\int_Sa(x,y)|\eta^1_t-\eta^2_t|(dy)]\label{cc311}\\
&&\qquad\le e^{-t}\sup_xc(x)E[\int_Sa(x,y)|\eta_0^1-\eta^2_0|(dy)
]\nonumber\\
&&\qquad\qquad +M\int_0^t\sup_x\,c(x)E[\int_Sa(x,y)|\eta^1_s-\eta^
2_s|(dy)]\,e^{-(t-s)}\,ds.\nonumber\end{eqnarray}

Multiply \reff{cc311} by $e^t$ and apply Gronwall's 
inequality to obtain the exponential rate of convergence.  
\end{proof}

\subsection{Spatial ergodicity} 

In this section, we take $S=\R^d$ and assume that $\lambda$ is 
translation invariant in the following sense.  For 
arbitrary $x,y\in\R^d$ and $B\in {\mathcal B}(\R^d)$, write 
\[T_xy=x+y~~~\mbox{\rm and}~~~T_xB=B+x=\{y+x;y\in B\}.\]
Then, $T_x$ induces a transformation $S_x$ on ${\mathcal N}(\R^d)$ through 
the equation 
\begin{equation}(S_x\eta )(B)=\eta (T_xB),~~~\eta\in {\mathcal N}(\R^
d),B\in {\mathcal B}(\R^d).\label{eqcolbsx}\end{equation}
Note that if $\eta =\sum\delta_{x_i}$, then $S_x\eta =\sum\delta_{
x_i-x}$.

\begin{definition} {\rm
We say that $\lambda$ is {\em translation invariant\/} if  
$\lambda (x+y,\eta )=\lambda (x,S_y\eta )$ for $x,y\in\R^d,\eta\in 
{\mathcal N}(\R^d)$.  }
\end{definition}

\begin{definition} {\rm
An ${\mathcal N}(\R^d)$-valued random variable $\eta$ is 
{\em translation invariant\/} if the distribution of $S_y\eta$ does not depend on 
$y$.  A probability distribution $\mu\in {\mathcal P}({\mathcal N}(\R^d))$ is 
{\em translation invariant\/} if  $\int f(\eta )\mu (d\eta )=\int 
f(S_y\eta )\mu (d\eta )$, 
for all $y\in {\Bbb R}^d$ and all bounded, measurable functions $
f$.  }
\end{definition} 

\begin{definition}\label{spergdef}
Let $\eta$ be a translation 
invariant, 
${\mathcal N}(\R^d)$-valued random variable.  
A measurable subset $G\subset {\mathcal N}(\R^d)$ is 
 {\em almost surely translation invariant\/} for $\eta$, if 
\[\one_G(\eta )=\one_G(S_x\eta )\quad a.s.\]
for every $x\in\R^d$. 
 $\eta$ is {\em spatially ergodic\/} if $P\{\eta\in G\}$ is $0$ 
or $1$ for each almost surely translation invariant $G\subset {\mathcal N}
(\R^d)$.
\end{definition}

Similarly, for $x\in\R^d$, we define $S_xN$ so that the spatial 
coordinate of each point is shifted by $-x$.  Almost sure 
translation invariance  of a set $G\subset {\mathcal N}(\R^d\times [0
,\infty )^3)$ and 
spatial ergodicity are defined analogously to Definition 
\ref{spergdef}.  Spatial ergodicity for $N$ follows from its 
independence properties.

\begin{lemma}\label{trinv}
Suppose $\lambda$ is translation invariant.  If $\eta_0$ is 
translation invariant and spatially ergodic and the 
solution of (\ref{stoch5}) is unique, then for each $t>0$, 
$\eta_t$ is translation invariant and spatially ergodic.
\end{lemma}

\begin{proof}\  $\{S_y\eta_t,t\geq 0\}$ is the solution of (\ref{stoch5})  
with $\eta_0$ replaced by $S_y\eta_0$ and $N$ replaced by $S_yN$.  By 
uniqueness,  $S_y\eta_t$ must have the same distribution as $\eta_
t$ 
giving the stationarity.  Also, by uniqueness, for 
measurable $G\subset {\mathcal N}(\R^d)$ there exists a measurable 
$\hat {G}\subset {\mathcal N}(\R^d)\times {\mathcal N}(\R^d\times [0,\infty 
)^3)$ such that 
\[\one_{\{S_x\eta_t\in G\}}=\one_{\{(S_x\eta_0,S_xN)\in\hat {G}\}}
\quad a.s.\]
for all $x\in\R^d$.  Consequently, spatial ergodicity for $\eta_t$ 
follows from the spatial ergodicity
of  $(\eta_0,N)$.
\end{proof}

\begin{remark}
{\rm If  $\eta$ is temporally ergodic and $\pi$ is the unique 
stationary distribution, then it must be translation 
invariant since $\{\eta_t\}$ stationary (in time) implies $\{S_x\eta_
t\}$ is 
stationary.
}
\end{remark}

\begin{lemma}\label{sperg}
Suppose that $\lambda$ is translation invariant and
satisfies (\ref{biglam}) and Condition 
\ref{contcnd}. 
 Then for each $t$, $\eta^{1,\infty}_t\equiv\lim_{n\rightarrow\infty}
\eta_t^{1,\infty ,n}$ and 
$\eta^{2,\infty}_t\equiv\lim_{n\rightarrow\infty}\eta_t^{2,\infty 
,n}$ are spatially ergodic.  

\end{lemma}

\begin{proof}
As in (\ref{map2}), $\eta_t^{1,\infty}=H^1(R_tN)$
 can be written as a 
deterministic transformation $F(t,N)$ of $N$ and that 
$S_y\eta_t^{1,\infty}=F(t,S_yN)$.  The spatial ergodicity of $\eta_
t^{1,\infty}$ then 
follows from the spatial ergodicity of $N$.\hfill \end{proof}

\begin{corollary}
If in 
addition to the conditions of Lemma \ref{sperg}, 
$\lambda$ satisfies the conditions of Lemma 
\ref{uniquesd}, then the unique stationary distribution is 
spatially ergodic.  In particular, if $\lambda$ satisfies the 
conditions of Theorem \ref{thmeu}\ with $M<1$, then the 
unique stationary distribution is spatially ergodic.
\end{corollary}

\begin{corollary}
If in addition to the conditions of Lemma \ref{sperg}, $\lambda$ is 
nondecreasing, then the minimal and maximal stationary 
distributions are spatially ergodic.
\end{corollary}

\vskip5mm 

\noindent {\bf Acknowledgments} This project was conducted 
during several visits of NLG to the CMS - University of Wisconsin.
This project was partially supported by CNPq Grant 301054/93-2 and
FAPESP 1995/4996-3 (NLG) and DMS 02-05034 and DMS
05-03983 (TGK). This material is based upon work supported by, or
in part by, the U.  S.  Army Research Laboratory and the U.  S.  Army
Research Office under contract, grant number DAAD19-01-1-0502.  and by
NSF Grant DMS 02-05034. \\

\end{document}